\documentclass[12pt]{amsart}

\usepackage{amscd,latexsym,amsthm,amsfonts,amssymb,amsmath,amsxtra}
\usepackage[mathscr]{eucal}
\usepackage{tkz-euclide}
\usepackage{pst-plot}
\usepackage{hyperref}
\renewcommand\theequation{\thesection.\arabic{equation}}

\newcommand{\BA}{{\mathbb {A}}}

\newcommand{\CA}{{\mathcal {A}}}

\newcommand{\CE}{{\mathcal {E}}}

\newcommand{\I}{{\mathrm{I}}}

\newcommand{\bs}{\backslash}

\newtheorem*{theorem*}{Theorem}
\newtheorem{lemma}{Lemma}[section]
\newtheorem{proposition}[lemma]{Proposition}

\newtheorem{theorem}[lemma]{Theorem}

\newtheorem{corollary}[lemma]{Corollary}

\newtheorem*{conjecture*}{Conjecture}

\makeatletter

\newcommand{\Rmnum}[1]{\expandafter\@slowromancap\romannumeral #1@}
\makeatother

\begin{document}
\renewcommand{\theequation}{\arabic{equation}}
\numberwithin{equation}{section}

\title[Global Converse Theorem]{A remark on a converse Theorem of Cogdell and Piatetski-Shapiro}

\author{Herv\'e Jacquet}
\address{Department of Mathematics\\
Columbia University\\
Rm 615, MC 4408 2990 Broadway\\
New York, NY 10027 USA}
\email{hj@math.columbia.edu}

\author{Baiying Liu}
\address{Department of Mathematics\\
Purdue University\\
150 N. University St\\
West Lafayette, IN, 47907 USA}
\email{liu2053@purdue.edu}

\begin{abstract}
In this paper, we reprove a global converse theorem of Cogdell and Piatetski-Shapiro using purely global methods.  
\end{abstract}

\date{\today}
\thanks{The second mentioned author was supported in part by NSF Grant DMS-1620329 and in part by a start-up fund from the Department of Mathematics of Purdue University.}
\maketitle


\section{Introduction}
Let $F$ be a number field or a function field. Denote by $\mathbb{A}$ the ring of adeles of $F$ and by $\psi$ a non-trivial additive character of $F \bs \mathbb{A}$.  Let $n\geq 4$. Let $\pi$ be an irreducible generic representation of $GL_n(\mathbb{A})$. We assume that the central character $\omega_\pi$ of $\pi$ is automorphic (condition $\CA(n,0)$). We also assume that if $\tau$ is a cuspidal automorphic representation of $GL_m(\mathbb{A})$ the complete $L-$function $L(s,\pi\times  \tau)$ converges for $\mathrm{Re}  s$ large enough. We denote by $\mathcal{A}(n,m)$ the condition that, for every such $\tau$, the $L-$function $L(s,\pi \times \tau)$ has the standard analytic properties (is nice in the terminology of Cogdell and Piatetski-Shapiro \cite{CPS94, CPS96, CPS99, Cog02}, see p. 4 for details). 

Following them, for every $\xi$ in the space $V_{\pi}$ of $\pi$, we let $W_\xi$ be the corresponding element of the Whittaker model $\mathcal{W}(\pi,\psi)$ of $\pi$.  We denote by $U_n$ the group of upper triangular matrices in $GL_n$ with unit diagonal. We set
\[ U_\xi (g) = \sum _{\gamma \in U_{n-1}(F)\backslash GL_{n-1}(F)} W_\xi \left[\left(\begin{array}{cc} \gamma & 0 \\ 0 & 1 \end{array}\right) g\right]\,,\]
\[ V_\xi(g) =  \sum _{\gamma \in U_{n-1}(F)\backslash GL_{n-1}(F)} W_\xi\left[ \left(\begin{array}{cc} 1 & 0 \\ 0 & \gamma \end{array}\right) g\right]\,.\]
If $\pi$ is automorphic cuspidal then $U_\xi=V_\xi$ for all $\xi \in V_{\pi}$. Conversely, if $U_\xi=V_\xi$ for all $\xi \in V_\pi$ or, what amounts to the same,  $U_\xi(I_n)= V_\xi(I_n)$ for all $\xi \in V_{\pi}$, then $\pi$ is automorphic. 
 
Let $Z_{U_n}$ be the center of the group $U_n$. Cogdell and Piatetski-Shapiro (\cite{PS76}, \cite{CPS96}, \cite{CPS99}) prove that the conditions $\mathcal{A}(n,m)$ with $0\leq m\leq n-2$ imply that
\[ \int _{Z_{U_n}(F) \bs Z_{U_n}(\mathbb{A})} (U_\xi - V_\xi) (z) \theta(z) d z =0\]
for all non-trivial characters $\theta$ of $Z_{U_n}(F)\bs Z_{U_n}(\mathbb{A})$. 
They do not have the same relation for the trivial character  which would then imply that $U_\xi(I_n)=V_\xi(I_n)$ for all $\xi \in V_{\pi}$. Nonetheless, they prove that $\pi$ is automorphic by using an ingenious local construction.

Our goal in this paper is to prove that conditions $\mathcal{A}(n,m)$, $0\leq m\leq n-3$, imply that
\[ \int _{Z_{U_n}(F) \bs Z_{U_n}(\mathbb{A})}  (U_\xi - V_\xi)( z) d z=0 \,,\, \forall \xi \in V_{\pi} \,.\]
This proves directly that the conditions $\mathcal{A}(n,m)$ with $0\leq m\leq n-2$ imply $U_\xi(I_n)=V_\xi(I_n)$ for all $\xi \in V_{\pi}$ and, in turn, imply $\pi$ is automorphic (and cuspidal). 

While our result is not needed, it gives a purely global proof of the Theorem of Cogdell and Piatetski-Shapiro. It is also germane to the conjecture that the conditions $\mathcal{A}(n,m)$ with $0\leq m\leq [\frac{n}{2}]$ imply that $\pi$ is automorphic (and cuspidal). Of course, the conjecture is true for $n=2,3, 4$ (see \cite{JL70}, \cite{JPSS79}, \cite{PS76} and \cite{CPS96}). 

The material is arranged as follows. In the next section, for the convenience of the reader, we review the work of Cogdell and Piatetski-Shapiro and state our result.  In section $3$ we provide preliminary material of an elementary nature. In section 4 we prove our result.

\textbf{Acknowledgements.} 
The second mentioned author would like to thank James Cogdell for helpful conversation on their previous results when he was visiting Ohio State University. 
The authors also would like to thank him for helpful comments and suggestions on an earlier version of the paper. 
This material is based upon work supported by the National Science Foundation under agreement No. DMS-1128155. Any opinions, findings and conclusions or recommendations expressed in this material are those of the authors and do not necessarily reflect the views of the National Science Foundation. We also would like to thank the referee for a very careful reading of the paper and for many helpful comments and suggestions.

\section{Preliminaries and the main result}
In $G_n=GL_n$ we let $U_n$ be the subgroup of upper triangular matrices with unit diagonal. We let $A_n$ be the group of diagonal matrices and $Z_n$ the center of $G_n$.
We define a character $\psi_{U_n}$ of $U_n(\mathbb{A})$ which is trivial on $U_n(F)$ by
\[ \psi_{U_n}( u) = \psi( u_{1,2}+ u_{2,3}+\cdots + u_{n-1,n}) \,.\]
We let $\pi$ be an irreducible generic representation of $G_n(\mathbb{A})$. As usual, this means that $\pi$ is a restricted tensor product of local irreducible representations $\pi_v$. For a finite place $v$, $\pi_v$ is an irreducible admissible representation of $G_n(F_v)$ on a complex vector space $V_v$. We assume that $\pi_v$ is generic, that is, there is a non-zero linear form
\[ \lambda_v: V_v \rightarrow  \mathbb{C} \]
such that 
\[ \lambda_v( \pi_v(u) e ) = \psi_{U_n,v}(u) \lambda_v(e) \]
for all vectors $e\in V_v$ and all $u\in U_n(F_v)$. We denote by $\mathcal{W}( \pi_v,\psi_{U_n,v})$ the space of functions
\[ g \mapsto \lambda _v( \pi_v(g) e ) \,,\, e\in V_v\]
on $G_n(F_v)$. It is the Whittaker model of $\pi_v$ noted $\mathcal{W}(\pi_v, \psi_{U_n,v})$.
For all finite $v$ not in a finite set $S$, the space contains a unique vector $W_{v,0}$ fixed under $G_n(\mathcal{O}_v)$ and taking the value $1$ at $I_n$. 
The representation $\pi_v$ is then determined by its Langlands semi-simple conjugacy class $A_v\in G_n(\mathbb{C})$. We assume that there is an integer $m\geq 0$ such that for all finite $v\not\in S$, any eigenvalue $\alpha$ of $A_v$ verifies $q_v^{-m} \leq |\alpha | \leq q_v^{m}$. 

For an infinite place $v$, the representation $\pi_v$ is really an irreducible admissible Harish-Chandra module. We denote by $(\pi_v,V_v)$ its canonical completion of slow growth in the sense of Casselman and Wallach. We assume that there is a non-zero continuous linear form 
\[ \lambda_v: V_v \rightarrow  \mathbb{C} \]
satisfying the same condition as before. We also define $\mathcal{W}(\pi_v,\psi_{U_n,v})$ as before.

Finally, let $\infty$ be the set of infinite places of $F$ and
\[ G_{n, \infty}= \prod _{v \in \infty} G_n(F_v) \]
We let $(\pi_\infty, V_\infty)$ be the topological tensor product of the representations $(\pi_v, V_v)$, $v\in \infty$. Let $\lambda$ be the tensor product of the linear forms $\lambda_v$, $v\in \infty $. We can define the space $\mathcal{W}(\pi_\infty,\psi_{U_n,\infty})$. 

We denote by $V$ the restricted tensor product $\bigotimes'_{v} V_v$ and $\pi$ the natural representation of $G_n(\mathbb{A})$ on $V_{\pi}=V$. 
The Whittaker model $\mathcal{W}(\pi,\psi_{U_n})$ of $\pi$ is the space spanned by the functions 
\[  W_\infty \prod _{v\not\in \infty}W_v \]
with $W_\infty\in \mathcal{W}(\pi_\infty,\psi_{U_n,\infty})$, $W_v\in \mathcal{W}(\pi_v,\psi_{U_n,v})$ and $W_v =W_{v,0}$ for almost all $v\not\in S$.  
For every $\xi\in V_{\pi}$ we denote by $W_\xi$ the corresponding element of $\mathcal{W}(\pi,\psi_{U_n})$. 

We assume that the central character of $\pi$ is automorphic. It is convenient to refer to this condition as condition $\mathcal{A}(n,0)$.
In view of our assumptions, for any cuspidal automorphic representation $\tau$ of $G_m(\mathbb{A})$, $1 \leq m\leq n-1$, the complete $L-$function $L(s, \pi \times \tau)$ is defined by a convergent product for $\mathrm{Re} s$ large enough. {\it Condition $\mathcal{A}(n,m)$} is that, for any such $\tau$, the function $L(s,\pi\times \tau)$ extends to an entire function of $s$, bounded in vertical strips and satisfies the functional equation
\[ L(s, \pi \times \tau) = \epsilon(s, \pi \times \tau, \psi) L(1-s, \widetilde{\pi} \times \widetilde{\tau}) \,.\]

In $G_n$ let $Y_{n,m}$ be the unipotent radical of the standard parabolic subgroup of type $(m+1, 1, 1,\ldots,1)$. For instance:
\[ Y_{3,1} = \left\{ \left( \begin{array}{ccc}  1 & 0 & \bullet \\ 0 & 1 & \bullet \\ 0 &0 & 1 \end{array}\right) \right\}\,,\,
Y_{5,2} = \left\{ \left(\begin{array}{ccccc} 
1 & 0& 0 &  \bullet &\bullet \\
0&1 & 0 & \bullet &\bullet \\
0&0&1 & \bullet  &\bullet \\
0&0&0 & 1 & \bullet \\
0&0&0&0& 1
\end{array}\right) \right\}\,. \]
If a function $\phi$ on $G_n(\mathbb{A})$ is invariant on the left under $Y_{n,m}(F)$ we set 
\[ \mathbb{P}^n_m (\phi)( g) = \int_{Y_{n,m}(F)\backslash Y_{n,m}(\mathbb{A})} \phi (y g) \overline{\psi}_{U_n}(y) d y  \,.\]
Here $ d y$ is the Haar measure on $ Y_{n,m}(\mathbb{A})$ normalized by the condition that the quotient $Y_{n,m}(F)\backslash Y_{n,m}(\mathbb{A})$ has measure $1$. Our notation differs slightly from the notations of Cogdell and Piatetski-Shapiro (\cite{Cog02}). Here $\mathbb{P}^n_m (\phi)$ is a function on $G_n(\BA)$, while in \cite{Cog02}, it is a function on $P_{m+1}(\BA)$, the mirabolic subgroup of $G_{m+1}(\BA)$, embedded in $G_n(\BA)$.
Note that $Y_{n,n-1}=\{ I_n\}$ and $\mathbb{P}^n_{n-1}$ is the identity.

Suppose $\pi$ is as above. For each $\xi$ in the space $V_{\pi}$ of $\pi$ we set
\[ U_\xi (g) = \sum _{\gamma \in U_n(F) \backslash P_n(F)} W_\xi(\gamma g),\]
where $P_n$ is the subgroup of matrices of $G_n$ whose last row has the form
\[ (0,0,\ldots,0,1) \,.\]
This is also
 \[U_\xi (g) =\sum _{\gamma \in U_{n-1}(F) \backslash G_{n-1}(F)} W_\xi\left[ \left(\begin{array}{cc} \gamma& 0 \\ 0 & 1 \end{array}\right) g\right]\,.\]
Then we have the following result.

\begin{lemma}
With the previous notations, 
\[ \mathbb{P}^n_m(U_\xi)( g) = \sum _{\gamma \in U_{m+1}(F)\backslash P_{m+1}(F)} W_{\xi} \left[ \left(\begin{array}{cc} \gamma &0 \\ 0 & I_{n-m-1} \end{array}\right) g\right]\,,\] 
or equivalently
\[\mathbb{P}^n_m(U_\xi)( g) =  \sum _{\gamma \in U_{m}(F)\backslash G_{m}(F)} W_{\xi} \left[ \left(\begin{array}{cc} \gamma &0 \\ 0 & I_{n-m} \end{array}\right) g\right] \,.\]
\end{lemma}

Likewise, let $R_n$ be the subgroup of matrices of $G_n$ whose first column is
\[ \left(\begin{array}{c} 1 \\ 0 \\0 \\\vdots \\0 \end{array}\right) \,.\]
We can consider the function
\[ V_\xi (g) = \sum _{\xi \in U_n(F)\backslash R_n(F)} W_\xi (\gamma g) \,.\]  
Let $\omega_n$ be the permutation matrix defined by
\[ \omega_1=1\,,\, \omega_n = \left(\begin{array}{cc} 0 & \omega_{n-1} \\ 1 & 0 \end{array}\right)\,.\]
Then
\[ R_n =\omega_n  {}^ t P_n^{-1} \omega_n \,,\, U_n =\omega_n {}^t U_n^{-1} \omega_n \,.\]
Moreover the automorphism $u\mapsto \omega_n {}^t u^{-1} \omega_n$ changes $\psi_{U_n}$ into $\overline{\psi}_{U_n}$.

If $\pi$ is automorphic cuspidal then for the cusp form $\phi_\xi$ corresponding to $\xi \in V_{\pi}$
we have 
\[ W_\xi(g) = \int_{U_n(F) \bs U_n(\BA)} \phi_\xi( u g) \overline{\psi}_{U_n}(u) d u \]
and 
\[ \phi_\xi(g)= U_\xi(g)\,. \]
By the previous observation relative to $R_n$, we also have
\[  \phi_\xi(g)=  V_\xi(g)\,. \]
Thus we have
\[ U_\xi = V_\xi \,,\, \forall \xi  \in V_{\pi} \,.\]

Conversely if $\pi$ is given and
\[  U_\xi = V_\xi \,,\, \forall \xi  \in V_{\pi} \,,\]
then $U_\xi$ is invariant on the left under $Z_n(F), P_n(F), R_n(F)$. Since these groups generate $G_n(F)$, for every $\xi \in V_{\pi}$ the function $U_\xi$ is invariant on the left under $G_n(F)$ and hence $\pi$ is automorphic.  

In general, $U_\xi$ and $V_\xi$ are invariant on the left under $P_n(F)\cap R_n(F)$ and $A_n(F)$. 
In other words, they are invariant under $S_n(F)$ where $S_n$ is the standard parabolic subgroup of type $(1,n-2,1)$. The notations here differ slightly from those of Cogdell and Piatetski-Shapiro (\cite{CPS99}). 
Moreover, we have, for all $g, h \in G_(\BA)$, 
\[W_{\pi(h)\xi} (g) = W_\xi( g h )\]
and similar formulae for $U_{\xi}$ and $V_{\xi}$. As a consequence, if an identity involving $W_{\xi}$, $U_{\xi}$ or $V_{\xi}$ is true for all $\xi \in V_{\pi}$, the identity obtained by translating the function on the right by an arbitrary element of $G_n(\BA)$ is also true for all $\xi \in V_{\pi}$, and conversely. 
For instance, the relation $U_\xi(I_n) = V_\xi(I_n)$ for all $\xi \in V_{\pi}$ is equivalent to the relation $U_\xi(g) =V_\xi(g)$ for all $\xi \in V_{\pi}$ and for all $g \in G_n(\BA)$. We appeal repeatedly to this principle.

Following Cogdell and Piatetski-Shapiro (\cite[Section 5.2]{Cog02}), for $1\leq m\leq n-2$, we define
\[ \alpha_m = \left( \begin{array}{ccc}0 & 1 & 0 \\I_m&0&0 \\ 0&0 &I_{n-m-1}\end{array}\right) \,.\] 
We note that $\alpha_m \in P_n$.
We also define
\[ V_\xi^m(g) = V_\xi( \alpha_m g) \,, \forall g \in G_n(\BA)\,.\]
Thus $V_\xi^m$ is invariant under $Q_m = \alpha_m^{-1}R_n \alpha_m$.  This is the subgroup of matrices of $G_n$ whose $(m+1)$-th column has the form
\[ \left(\begin{array}{c} 0\\ \vdots \\ 0\\ 1 \\0 \\\vdots\\ 0\end{array}\right)\,, \]
with $1$ in the $(m+1)$-th row.  Note that $Q_m$ contains the group $Y_{n,m}$. Thus we may consider $ \mathbb{P}^n_m(V_\xi^m)$. 

\begin{theorem}[Cogdell, Piatetski-Shapiro, Section 5 of \cite{Cog02}, \cite{CPS99}]
Suppose conditions $\mathcal{A}(n,j)$, $0\leq j\leq m$, are satisfied. Then, for all $\xi \in V_{\pi}$, 
$$\mathbb{P}^n_m(U_\xi)= \mathbb{P}^n_m(V_\xi^m).$$
\end{theorem}

We denote by $\CE(n,m)$ the condition that
\[ \mathbb{P}^n_m(U_\xi)= \mathbb{P}^n_m(V_\xi^m)\,,\, \forall \xi \in V_\pi \,.\]
This condition can be simplified. Indeed, we have the following result. 

\begin{proposition}[Cogdell, Piatetski-Shapiro, Section 5 of \cite{Cog02}, \cite{CPS99}]\label{equivalentprop}
Let $k, 1 \leq k\leq n-m$, be an integer. The condition 
$\CE(n,m)$
is equivalent to the condition 
\begin{align}\label{condition1}
\begin{split}
& \ \int_{U_{n-m}(F)\backslash U_{n-m}(\mathbb{A})} 
\int
U_\xi  \left[ \left(\begin{array}{cc} I_{m} & x \\ 0 & u\end{array}\right) g \right]\overline{\psi}_{U_{n-m}}(u) d x d u\\
= & \ \int_{U_{n-m}(F)\backslash U_{n-m}(\mathbb{A})} \int V_\xi ^m \left[ \left(\begin{array}{cc} I_{m} & x \\ 0 & u\end{array}\right) g \right]\overline{\psi}_{U_{n-m}}(u) d x d u\,,
\end{split}
\end{align}
for all $\xi \in V_{\pi}$ and for all $g \in G_n(\BA)$,
where $x \in M_{m \times (n-m)}(F) \bs M_{m \times (n-m)}(\BA)$ with zero first $k$ columns. 
\end{proposition}

\textsc{Proof:} For the convenience of the reader, we review the proof. For $k=1$ our conclusion is just the hypothesis. Thus we may assume  our assertion true for $k$, $1 \leq k \leq n-m-1$ and prove it for $k+1$. In the integral we write
$ x = (0,0,\ldots,0, x_{k+1},x_{k+2},\ldots x_{n-m})$ where the $x_i$ are column vectors of length $m$. We also introduce
\[ \gamma_\beta = \left(\begin{array}{ccc}
I_m &0&0 \\ X_\beta & I_k &0 \\0 &0 & I_{n-m-k} 
\end{array}\right)\,,\, X_\beta= \overbrace{\left. \left(\begin{array}{c} 0\\0\\ \vdots \\0 \\ \beta\end{array}\right)\right\}}^m k \,,\, \beta\in F^m\,.
\] 
Since $\gamma_\beta$ is in $P_n(F)\cap Q_m(F)$, in the identity, we can conjugate the matrices by $\gamma_\beta$. We note that
\[ \gamma_\beta \left( \begin{array}{cc} I_m & x \\0 & u\end{array}\right)\gamma_\beta^{-1}=
 \left( \begin{array}{cc} I_m & x \\0 & u'\end{array}\right)\]
where
\[ u'_{m+k, m+k+1} =  u_{m+k, m+k+1}+\beta x_{k+1} \,.\]
Hence the equality becomes
\begin{align*}
& \ \int_{U_{n-m}(F)\backslash U_{n-m}(\mathbb{A})} 
\int
U_\xi  \left[ \left(\begin{array}{cc} I_{m} & x \\ 0 & u\end{array}\right) g \right]\overline{\psi}_{U_{n-m}}(u) \psi( -\beta x_{k+1})d x d u\\
= & \ \int_{U_{n-m}(F)\backslash U_{n-m}(\mathbb{A})} \int V_\xi ^m \left[ \left(\begin{array}{cc} I_{m} & x \\ 0 & u\end{array}\right) g \right]\overline{\psi}_{U_{n-m}}(u)\psi( -\beta x_{k+1}) d x d u\,,
\end{align*}
where $x \in M_{m\times (n-m)}(F) \backslash M_{m\times (n-m)}(\mathbb{A})$, with zero first $k$ columns. Summing over all $\beta \in F^m$ and applying the theory of Fourier series, we get
\begin{align*}
& \ \int_{U_{n-m}(F)\backslash U_{n-m}(\mathbb{A})} 
\int
U_\xi  \left[ \left(\begin{array}{cc} I_{m} & x \\ 0 & u\end{array}\right) g \right]\overline{\psi}_{U_{n-m}}(u) d x d u\\
= & \ \int_{U_{n-m}(F)\backslash U_{n-m}(\mathbb{A})} \int V_\xi ^m \left[ \left(\begin{array}{cc} I_{m} & x \\ 0 & u\end{array}\right) g \right]\overline{\psi}_{U_{n-m}}(u) d x d u\,,
\end{align*}
where $x \in M_{m\times (n-m)}(F) \backslash M_{m\times (n-m)}(\mathbb{A})$, with zero first $k+1$ columns.

The other direction is obvious, since the if the condition \eqref{condition1} holds for $k+1$, then integrating both sides with respect to 
$$x \in M_{m \times (n-m)}(F) \bs M_{m \times (n-m)}(\BA)$$
with only $(k+1)$-th column being possibly nonzero, we obtain the condition \eqref{condition1} for $k$, which is equivalent to the condition $\CE(n,n-m)$ by induction assumption.  
This concludes the proof of the proposition. \qed

Since $\alpha_m\in P_n$ and $U_\xi$ is $P_n(F)$ invariant on the left, we can apply the condition \eqref{condition1} with $g$ replaced by $\alpha_m^{-1}$. Thus the condition \eqref{condition1} can be written also as
\begin{align}\label{condition2}
\begin{split}
& \ \int_{U_{n-m}(F)\backslash U_{n-m}(\mathbb{A})} 
\int
U_\xi  \left[ \alpha_m\left(\begin{array}{cc} I_{m} & x \\ 0 & u\end{array}\right) \alpha_m^{-1} \right]\overline{\psi}_{U_{n-m}}(u) d x d u\\
= & \ \int_{U_{n-m}(F)\backslash U_{n-m}(\mathbb{A})} \int V_\xi \left[ \alpha_m\left(\begin{array}{cc} I_{m} & x \\ 0 & u\end{array}\right) \alpha_m^{-1} \right]\overline{\psi}_{U_{n-m}}(u) d x d u\,,
\end{split}
\end{align}
for all $\xi \in V_{\pi}$,
where $x \in M_{m \times (n-m)}(F) \bs M_{m \times (n-m)}(\BA)$ with zero first $k$ columns. 
In this paper, we are mainly interested in the case $k=n-m$ of condition \eqref{condition2}. We make it explicit.

Taking $k= n-m$ and $m=n-2$, then condition \eqref{condition2} leads to the condition of Cogdell and Piatetski-Shapiro
\[ \int_{F\backslash \mathbb{A}} (U_\xi - V_\xi) \left( \begin{array}{ccc}  1 & 0 & z \\0 & I_{n-2} & 0 \\ 0 & 0 & 1 \end{array}\right) \psi(-z) d z =0\]
for all $\xi \in V_{\pi}$. 
Since $U_\xi$ and $V_\xi$ are invariant on the left under $A_n(F)$ and this relation is true for any right translate of $U_\xi$ and $V_\xi$ we can conjugate by an element of $A_n(F)$ and obtain the condition
\[ \int_{F\backslash \mathbb{A}} (U_\xi - V_\xi) \left( \begin{array}{ccc}  1 & 0 & z \\0 & I_{n-2} & 0 \\ 0 & 0 & 1 \end{array}\right) \psi(-\alpha z) d z =0\,,\]
for all $\xi \in V_{\pi}$ and all $\alpha\in F^\times$. 
 
Now let us take $k=n-m$ and $m= n-3$.  Let $e_0$ and $f_0$ be respectively the following row and column of size $n-2$:
\[ e_0 = (0,0,\ldots,0,1) \,,\, f_0 =\left(\begin{array}{c} 0 \\0 \\ \vdots \\ 0\\1 \end{array}\right) \,.\]
Condition \eqref{condition2} reads
\[ \int_{(F\backslash \mathbb{A})^2} \int_{F\backslash \mathbb{A}} (U_\xi - V_\xi) \left( \begin{array}{ccc}  1 & x e_0 & z \\0 & I_{n-2} &  y f_0 \\ 0 & 0 & 1 \end{array}\right) \psi(-x-y) d z d x d y =0 \,,\]
for all $\xi \in V_{\pi}$. Abusing notation, we write the above integral as
\[ \int_{(F\backslash \mathbb{A})^3} (U_\xi - V_\xi) \left( \begin{array}{ccc}  1 & x e_0 & z \\0 & I_{n-2} &  y f_0 \\ 0 & 0 & 1 \end{array}\right) \psi(-x-y) d z d x d y =0 \,,\]
for all $\xi \in V_{\pi}$. 
For instance, for $n=4$ the condition reads
\[ \int_{(F\backslash \mathbb{A})^3} (U_\xi - V_\xi) \left( \begin{array}{cccc}  1 & 0& x  & z \\0 &1 & 0 & 0 \\ 0 & 0 & 1 & y \\ 0 & 0 & 0 & 1\end{array}\right) \psi(-x-y) d z d x d y =0\,,\]
for all $\xi \in V_{\pi}$. 
Again we can conjugate by an element of $A_n(F)$ to obtain 
\begin{equation}\label{equvan}
\int_{(F\backslash \mathbb{A})^3} (U_\xi - V_\xi) \left( \begin{array}{ccc}  1 & x e_0 & z \\0 & I_{n-2} &  y f_0 \\ 0 & 0 & 1 \end{array}\right) \psi(-\alpha x-\beta y) d z d x d y =0\,,
\end{equation} 
for all $\xi \in V_{\pi}, \alpha \in F^\times, \beta \in F^\times$.
Our own contribution is the following.

\begin{theorem} \label{MainResult}
Suppose $n>3$. Then condition $\mathcal{E}(n,n-3)$ is equivalent to the condition
 \[ \int_{F\backslash \mathbb{A}} (U_\xi - V_\xi) \left( \begin{array}{ccc}  1 & 0 & z \\0 & I_{n-2} & 0 \\ 0 & 0 & 1 \end{array}\right)  d z =0\]
for all $\xi \in V_{\pi}$.
\end{theorem}

Combining our Theorem $\ref{MainResult}$ with the results of Cogdell and Piatetski-Shapiro (\cite{Cog02}) we arrive at the following result.
\begin{corollary}
Suppose conditions $\mathcal{A}(n,i)$, $0\leq i\leq n-3$, are satisfied. Then for every $\xi \in V_{\pi}$
\[ \int_{F\backslash \mathbb{A}} (U_\xi - V_\xi) \left(\begin{array}{ccc} 1 & 0 & z \\ 0 & I_{n-2} & 0 \\ 0 & 0 & 1 \end{array} \right) d z =0 \, . \]
\end{corollary}

\section{Separable functions}
Let $U$ and $V$ be vector spaces over $F$. A function
\[ \Phi: U(F)\backslash U(\mathbb{A}) \times V(F)\backslash V(\mathbb{A}) \rightarrow \mathbb{C} \]
is said to be (additively) \textit{separable}, if there exist two functions $\Phi_1:U(F)\backslash U(\mathbb{A})\rightarrow \mathbb{C} $ and 
$\Phi_2:V(F)\backslash V(\mathbb{A})\rightarrow \mathbb{C} $ such that, for all $(u,v)\in  U(\mathbb{A})\times V(\mathbb{A})$, 
\[ \Phi(u,v) = \Phi_1(u) +\Phi_2(v) \,.\]
It amounts to the same to demand that, for all $(u,v)$,
\[ \Phi(u,v) =\Phi(u,0) + \Phi(0,v) -\Phi(0,0) \,.\]
In what follows, if $\Phi: U(F)\backslash U(\mathbb{A}) \times V(F)\backslash V(\mathbb{A}) \rightarrow \mathbb{C}$ is a function, an integral
\[ \int \int \Phi(u,v) d u d v \]
means that the integral is over the product $U(F)\backslash U(\mathbb{A}) \times V(F)\backslash V(\mathbb{A})$ and the Haar measure $ d u$ (resp. $d v$) is normalized by demanding that the quotients have volume $1$. This convention remains in force for the rest of this paper.

\begin{lemma}\label{SeparableLemma}
Suppose that
\[ \Phi :  F\backslash \mathbb{A} \times   F\backslash \mathbb{A} \rightarrow \mathbb{C} \]
is a smooth function such that, for all $\alpha\in F^\times$, $\beta \in F^\times$,
\[ \int \int \Phi(u,v) \psi(\alpha u + \beta v) d u d v =0 \,.\]
Then $\Phi$ is separable. 
\end{lemma}

\textsc{Proof:} We write the Fourier expansion of $\Phi$,
\begin{align*}
& \ \Phi(x,y) \\
= & \ \sum _{\alpha\in F, \beta \in F} \psi(\alpha x + \beta y) \left( \int \int \Phi(u,v) \psi(-\alpha u - \beta v) d u d v \right) \,.
\end{align*}
In view of the assumptions, we have
\begin{align*}
& \ \Phi(x,y) \\
= & \ \int \int \Phi(u,v) d u  d v\\
+ & \ \sum_{\alpha \in F^\times} \psi(\alpha x) \left( \int \int \Phi(u,v) \psi(-\alpha u ) d u d v \right)\\
+ & \ \sum_{\beta \in F^\times} \psi(\beta y) \left( \int \int  \Phi(u,v) \psi(-\beta v) d u d v \right)\,.
\end{align*}
The function on the right hand side of the equation is indeed separable. \qed

\begin{proposition}\label{Separable}
Suppose $U$ and $V$ are finite dimensional spaces over $F$ in duality by the bi-linear form $(u,v)\mapsto \langle u, v\rangle$. Suppose 
\[ \Phi: U(F)\backslash U(\mathbb{A}) \times V(F)\backslash V(\mathbb{A}) \rightarrow \mathbb{C} \]
is a smooth function with the following property: for all pairs $(e,f)\in U(F)\times V(F)$ with $\langle e,f\rangle =1$ we have
\[\int_{(F \bs \BA)^2} \Phi (u+ x e, v+  y f) \psi(\alpha x+ \beta y) d x d y =0\]
for all $u\in U(\mathbb{A})$, all $v\in V(\mathbb{A})$, all $\alpha \in F^\times$, and all $\beta \in F^\times$. Then $\Phi$ is separable.
\end{proposition}
\textsc{Proof:} If $\mathrm{dim}(U) =\mathrm{dim }(V)=1$, our assertion follows from the previous lemma. Thus we may assume that  $\mathrm{dim}(U) =\mathrm{dim }(V)=n+1$, $n>0$ and our assertion is true for dimension $n$. 
Let $e\in U(F),f\in V(F)$ with $\langle e,f\rangle=1$. Let $U_1$ be the subspace of $U$ orthogonal to $f$ and $V_1$ the subspace of $V$ orthogonal to $e$. By Lemma \ref{SeparableLemma}, for $u\in U_1(\mathbb{A}), v\in V_1(\mathbb{A})$ we have
\[ \Phi[ u+ s e ,v + t f] = \Phi[u+s e,v] + \Phi [u,v+ t f] - \Phi[u,v] \,.\]
Each one of the functions
\[ (u,v)\mapsto \Phi[u+s e,v]\,,\,(u,v)\mapsto  \Phi [u,v+ t f]\,,\, (u,v)\mapsto   \Phi[u,v]\] satisfies the assumptions of the proposition. By the induction hypothesis, the right hand side is equal to
\[ \Phi[u+s e,0] + \Phi[ s e ,v] -\Phi[s e,0] + \Phi[ u,t f] + \Phi[ 0, v+ t f]- \Phi[0, t f] \]
\[ -\Phi[u,0] -\Phi[ 0,v] + \Phi[0,0] \,.\]
Thus it suffices to show that $( u,t) \mapsto \Phi[ u,t f]$ and $(s ,v) \mapsto \Phi[ s e ,v]$ are separable functions. Let us show this is the case for the first function. Let $e_1,e_2,\ldots ,e_n$ be a basis of $U_1(F)$. Write
\[ u = s_1 e_1 + s_2 e_2 +\cdots + s_n e _n \,.\]
Now $\langle e_1+ e ,f\rangle =1$. Thus,  
\[ \Phi[ s_1 e_1 + s_2 e_2+\cdots s_n e_n+ s_1 e, v+ t f ] \]
must be separable as a function of $(s_1,t)$. All the terms on the right hand side (with $s=s_1$) have this property, except possibly the term
\[  \Phi[ u,t f]\,. \]
Thus this term must have this property as well.  Hence
\[ \Phi[ s_1 e_1 + s_2 e_2 +\cdots + s_n e _n ,tf] \]
is a separable function of the pair $(s_1,t)$. Likewise it is a separable function of each pair $(s_j,t)$, $1 \leq j \leq n$. By the lemma below it is a separable function of $((s_1,s_2,\ldots,s_n),t)$, that is, $\Phi[u,tf]$ is a separable function of $(u,t)$. \qed

\begin{lemma}
Suppose
\[ \Phi( (s_1,s_2,\ldots ,s_n),t) \]
is a function with the property that for each index $j$ it is a separable function of $(s_j,t)$. Then it is a separable function of the pair
\[ ((s_1,s_2,\ldots, s_n),t) \,.\]
\end{lemma}

\textsc{Proof:}  Our assertion is obvious if $n=1$. So we may assume $n>1$ and our assertion true for $n-1$. We have, by separability in $(s_1,t)$,
\[ \Phi((s_1,s_2,\ldots s_n),t) = \]
\[  \Phi((s_1,s_2,\ldots s_n),0) + \Phi((0,s_2,\ldots s_n),t) - \Phi( (0,s_2,\ldots s_n), 0) \,.\]
By the induction hypothesis the term $\Phi((0,s_2,\ldots s_n),t)$ is a separable function of the pair $(( s_2,\ldots s_n),t)$. Thus the right hand side is a separable function of the pair $((s_1,s_2,\ldots s_n),t)$. \qed

Finally, we have a simple criterion to decide whether a separable function vanishes.

\begin{proposition}\label{Vanishing2}
Suppose $\Phi$ is a separable smooth function on 
$$U(F)\backslash U(\mathbb{A}) \times V(F)\backslash V(\mathbb{A})\,,$$
where $(U,V)$ is a pair of vector spaces over $F$ in duality by the bilinear form $(u,v)\mapsto \langle u,v\rangle$. Suppose that
\[ \int \int \Phi(u,v) \psi( \langle u,f\rangle + \langle e,v\rangle) d u  d v =0\,,\]
if $e\in U(F)$, $f\in V(F)$, and either $e=0$ or $f=0$ (or both $e=0$ and $f=0$). Then $\Phi=0$.
\end{proposition}
\textsc{Proof:}
By assumption,
\[ \Phi(u,v) = \Phi(u,0) + \Phi(0,v) -\Phi(0,0) \,.\]
Taking $e=0$ and $f\neq 0$ we get
\[\int \Phi(u,0) \psi(\langle u,f\rangle) d u + \left(\int \Phi(0,v) dv-\Phi(0,0)\right) \int \psi(\langle u,f \rangle)du=0 \,. \]
Since $f \neq 0$, there exists $u_0 \in U(F) \bs U(\BA)$, such that 
$\psi(\langle u_0,f \rangle) \neq 1$. By changing of variables, 
\[\int \psi(\langle u,f \rangle)du = \int \psi(\langle u+u_0,f \rangle)du= \psi(\langle u_0,f \rangle)\cdot \int \psi(\langle u,f \rangle)du \,,\]
which implies that $\int \psi(\langle u,f \rangle) du=0$. Hence we get
\[ \int \Phi(u,0) \psi(\langle u,f\rangle) d u =0 \,.\]
This shows that $u\mapsto \Phi(u,0)$ is a constant function. Likewise $v\mapsto \Phi(0,v)$ is a constant function. By the above formula from assumption, $(u,v) \mapsto \Phi(u,v)$ is a constant function. 
Now if we take $e=0$ and $f=0$ we get
\[\int \int \Phi(u,v) d u d v =0\,.\]
But since $\Phi$ is constant, this integral is just the constant value of $\Phi$. Hence $\Phi$ is $0$ as claimed. \qed

\section{Proof of Theorem \ref{MainResult}} 

It will be convenient to introduce for every $\xi\in V_\pi$ the function
\[ \Phi_\xi(u,v) = \int _{F\backslash \mathbb{A}}( U_\xi-V_\xi)\left(\begin{array}{ccc}
 1 & u & z \\ 0 & I_{n-2} & v \\ 0&0& 1 \end{array}\right) d z \,.\]
Here $u$ is a row of size $n-2$ and $v$ a column of size $n-2$. The scalar product of $u$ and $v$ is denoted by $\langle u, v\rangle$. Thus $\Phi_\xi$ is a smooth function on $(F\backslash \mathbb{A})^{n-2}  \times (F\backslash \mathbb{A})^{n-2}$. 

\begin{proposition} \label{Vanishing}
Let $A$ be a rational column of size $n-2$ and $B$ a rational row of size $n-2$. Suppose $A=0$ or $B=0$ (or both $A=0$ and $B=0$). Then, for every $\xi \in V_\pi$, the integral
\[ \int \int \Phi_\xi(u,v) \psi( \langle u,A\rangle+ \langle B, v \rangle ) d u d v \]
is $0$.
\end{proposition}

\textsc{Proof:} 
It suffices to prove the integral
\[\int \int\int  U_\xi\left(\begin{array}{ccc}
 1 & u & z \\ 0 & I_{n-2} & v \\ 0&0& 1 \end{array}\right) \psi( \langle u,A\rangle+
 \langle B,v\rangle) d u d v d z \]
has the properties described in the proposition. 

Indeed, the automorphism 
\[ g \mapsto \omega_n {}^t g^{-1} \omega_n \]
changes the function $U_\xi$ relative to the character $\psi$ into the function $V_\xi$
relative to the character $\psi^{-1}$ and leaves invariant the group over which we integrate.
Thus
the integral 
\[ \int\int \int V_\xi\left(\begin{array}{ccc}
 1 & u & z \\ 0 & I_{n-2} & v \\ 0&0& 1 \end{array}\right) \psi( uA+ Bv) d u d v d z \]
has the same properties and so does the integral of the difference $U_\xi-V_\xi$. Note that $\langle u, A \rangle = uA$, $\langle B, v \rangle = Bv$. 

Now consider the integral for $U_\xi$. If $B=0$ it suffices to show that, for every $\xi$,
\[  \int\int U_\xi\left(\begin{array}{ccc}
 1 & 0 & z \\ 0 & I_{n-2} & v \\ 0&0& 1 \end{array}\right) d v d z =0\,.\]
Indeed, as we have remarked, if this identity is true for all $\xi$, then the identity obtained by translating on the right by an element of $G_n$ is still true for every $\xi$. In particular, then,  for all $\xi$ and all $u$,
\[ \int \int  U_\xi\left[\left(\begin{array}{ccc}
 1 & 0 & z \\ 0 & I_{n-2} & v \\ 0&0& 1 \end{array}\right) 
\left(\begin{array}{ccc}
 1 & u & 0 \\ 0 & I_{n-2} & 0 \\ 0&0& 1 \end{array}\right)\right] 
 d v d z =0\,.\]
Integrating over $u$ with respect to the character $\psi(uA)$ we find
\[\int d u\int \int U_\xi\left[\left(\begin{array}{ccc}
 1 & 0 & z \\ 0 & I_{n-2} & v \\ 0&0& 1 \end{array}\right) 
\left(\begin{array}{ccc}
 1 & u & 0 \\ 0 & I_{n-2} & 0 \\ 0&0& 1 \end{array}\right)\right] 
 d v d z \psi(uA) du=0\,,\]
 or
 \[ \int \int \int  U_\xi\left(\begin{array}{ccc}
 1 & u & z \\ 0 & I_{n-2} & v \\ 0&0& 1 \end{array}\right)  \psi(uA)
d u d v d z =0\,,\]
as claimed.

Replacing $U_\xi$ by its definition, we find
\[  \int\int U_\xi\left(\begin{array}{ccc}
 1 & 0 & z \\ 0 & I_{n-2} & v \\ 0&0& 1 \end{array}\right) d v d z =\]
\[ \int  \int \sum _{\gamma  \in U_{n-1}(F) \backslash G_{n-1}(F)} W_\xi\left[ \left(\begin{array}{cc} \gamma & 0 \\ 0 & 1 \end{array} \right)\left(\begin{array}{ccc}
 1 & 0 & z \\ 0 & I_{n-2} & v \\ 0&0& 1 \end{array}\right) \right] d v d z\,.\]
Exchanging summation and integration we find
\[ \int \int \sum _{\gamma  \in U_{n-1}(F) \backslash G_{n-1}(F)} W_\xi\left[ \left(\begin{array}{cc} \gamma & 0 \\ 0 & 1 \end{array} \right)\left(\begin{array}{ccc}
 1 & 0 & z \\ 0 & I_{n-2} & v \\ 0&0& 1 \end{array}\right) \right] d v d z\,.\]
\[ =\sum _{\gamma  \in U_{n-1}(F) \backslash G_{n-1}(F)}\int \int   W_\xi\left[ \left(\begin{array}{cc} \gamma & 0 \\ 0 & 1 \end{array} \right)\left(\begin{array}{ccc}
 1 & 0 & z \\ 0 & I_{n-2} & v \\ 0&0& 1 \end{array}\right) \right] d v d z\,.\]
After a change of variables this becomes
\[  =\sum _{\gamma  \in U_{n-1}(F) \backslash G_{n-1}(F)}\int  \int  W_\xi\left[ \left(\begin{array}{ccc}
 1 & 0 & z \\ 0 & I_{n-2} & v \\ 0&0& 1 \end{array}\right)\left(\begin{array}{cc} \gamma & 0 \\ 0 & 1 \end{array}\right)\right] d v d z\]
which is $0$. 

Now suppose $B\neq 0$ (and thus $A=0$). To prove that the integral vanishes we may conjugate by a matrix
\[ \left(\begin{array}{ccc}
 1 & 0 & 0 \\ 0 & \gamma & 0 \\ 0&0& 1 \end{array}\right) \,,\, \gamma \in G_{n-2}(F)\,.\]
This amounts to replacing $B$ by $B \gamma^{-1}$. Thus we may assume $B= (0,0,\ldots,0,-1)$. Then the integral takes the form
\[ \int \mathbb{P}^n_{n-2}U_\xi \left( \begin{array}{ccc}
1 & u & 0 \\ 0 & I_{n-2} & 0 \\ 0&0& 1 \end{array}\right) d u \,.\]
Replacing $\mathbb{P}^n_{n-2}U_\xi $ by its expression in terms of $W_\xi$ we find
\[ \int \sum _{\gamma \in U_{n-2}(F) \backslash G_{n-2}(F)} W_\xi \left[  \left(\begin{array}{cc} \gamma & 0 \\ 0 & I_2 \end{array} \right) \left( \begin{array}{ccc}
1 & u & 0 \\ 0 & I_{n-2} & 0 \\ 0&0& 1 \end{array}\right)\right] d u \,.\]
Let us write this as
\[ \int \left( \int \sum_\gamma  W_\xi \left[  \left(\begin{array}{cc} \gamma & 0 \\ 0 & I_2 \end{array} \right) 
\left( \begin{array}{cccc}
1 & 0 & u'' & 0 \\ 0 & I_{n-3} & 0 &0\\ 0&0& 1 &0\\ 0 &0&0 &1  \end{array}\right)
\left( \begin{array}{ccc}
1 & u' & 0 \\ 0 & I_{n-3} & 0 \\ 0&0& I_2 \end{array}\right)\right] d u'' \right)d u' \,,\]
with $\gamma \in U_{n-2}(F) \backslash G_{n-2}(F)$. The inner integral is in fact a 
multiple of the integral
\[ \int \psi\left(\langle \gamma_{n-2}, \left(\begin{array}{c} u''\\ 0\end{array} \right)\rangle \right)d u' \,, \]
where $\gamma_{n-2}$ is the last row of $\gamma$. This integral is then $0$ unless the last row has the form
\[ (0,\bullet,\bullet,\ldots,\bullet) \,,\]
in which case the integral is $1$.
Thus our expression is
 \[ \int \sum_\gamma  W_\xi \left[  \left(\begin{array}{cc} \gamma & 0 \\ 0 & I_2 \end{array} \right) 
\left( \begin{array}{ccc}
1 & u' & 0 \\ 0 & I_{n-3} & 0 \\ 0&0& I_2 \end{array}\right)\right] d u' \,,\]
where $\gamma \in U_{n-2}(F) \backslash G_{n-2}(F)$ and the last row of $\gamma $ has the form $(0,\bullet,\bullet,\ldots,\bullet)$. 

Now let us write the Bruhat decomposition of $\gamma \in U_{n-2}(F) \backslash G_{n-2}(F)$:
\[ \gamma = w  \nu \alpha\]
with $w$ a permutation matrix, $\alpha$ a diagonal matrix and $\nu \in U_{n-2}(F)$. The last row of $w$ cannot have the form $(1,0,\ldots,0)$ otherwise the last row of $\gamma$ would have the form $(x,\bullet,\bullet,\ldots,\bullet)$, $x\neq 0$. 
Thus the last row of $w$ has the form $(0,\bullet,\bullet,\ldots,\bullet)$. Let us write down the contribution of such a $w$ to the above expression and show it is $0$. We introduce the abelian group
\[ X = \left\{ \left( \begin{array}{ccc}
1 & u' & 0 \\ 0 & I_{n-3} & 0 \\ 0&0& I_2 \end{array}\right)\right\} \,.\]
The contribution  of $w$ has the form 
\[ \int_{X(F)\backslash X(\mathbb{A})} \sum_\gamma  W_\xi \left[  \left(\begin{array}{cc} \gamma & 0 \\ 0 & I_2 \end{array} \right) x 
\right] d x \,,\]
where $\gamma = w  \nu \alpha$,   $\alpha \in A_{n-2}(F)$ and $\nu$ is in a set of representatives for the cosets $ w ^{-1} U_{n-2}(F)  \cap U_{n-2}(F) \backslash U_{n-2}(F)$. For a set of representatives, we will take the group 
\[S=  w ^{-1} \overline{U_{n-2}}(F) w \cap U_{n-2}(F) \,,\]
where as usual $\overline{U_{n-2}}$ is the subgroup opposite to $U_{n-2}$ (that is, its transpose). 

Now viewed as a subgroup of $G_{n-2}(F)$ the group $X(F)$ is the unipotent radical of the standard parabolic subgroup of type $(1,n-3)$ and $U_{n-2}(F)$ is contained in this parabolic subgroup. In particular, $X$ is normalized by $U_{n-2}(F)$ and by $S$. 

We keep in mind that $X$ is an abelian group. Let us write $X$ as the product $X_1 X_2$ 
where
\[ X_1(F) = w^{-1}U_{n-2}(F) w \cap X(F) \,,\, X_2(F) =  w^{-1}\overline{U_{n-2}}(F) w \cap X(F)
\,.\]
Thus $S$ contains $X_2(F)$ and is the product of $X_2(F)$ and the group $T$,
\[ T= S\cap \left\{ \left(\begin{array}{ccc} 
1&0&0\\
0 & \mu & 0 \\
0&0&I_2 \end{array}\right): \mu \in U_{n-3}(F) \right\}\,.\]
Moreover the group $S$ normalizes $X_2$. In particular it normalizes the groups
\[ X(F)= X_1(F)X_2(F)\,,\, X_2(\mathbb{A})\,, \]
hence also the closed subgroup
\[ X_1(F)X_2(\mathbb{A}) \,.\] 

Then our expression becomes
\[  \int \int \sum W_\xi\left[ \left(\begin{array}{cc} w & 0 \\ 0 & I_2 \end{array} \right)
\left(\begin{array}{cc}  \nu \alpha& 0 \\ 0 & I_2 \end{array} \right)  x_1 x_2  \right] d x_1 d x_2 \,.\]
Here $x_1$ and $x_2$ are integrated over $X_1(F)\backslash X_1(\mathbb{A})$ and $X_2(F)\backslash X_2(\mathbb{A})$ respectively. We sum for $\alpha \in A_{n-2}(F)$ and $\nu\in S$.
We can take the sum over $\alpha$ outside as follows:
\[ \sum _{\alpha}
\int \int \sum_\nu W_\xi\left[ \left(\begin{array}{cc} w & 0 \\ 0 & I_2 \end{array} \right)
\left(\begin{array}{cc} \nu & 0 \\ 0 & I_2 \end{array} \right)  x_1 x_2  
\left(\begin{array}{cc} \alpha  & 0 \\ 0 & I_2 \end{array} \right)\right] d x_1 d x_2 \,.\]
We now show that each term of the $\alpha$ sum is $0$ for all $\xi$. As usual, we may take $\alpha=I_{n-2}$. Now we write
\[ \nu = \tau \sigma\]
with $\sigma \in X_2(F)$ and $\tau$ in $T$. We combine the integration over $x_2$ and the sum over $\sigma$ to obtain an integral over $X_2(\mathbb{A})$. We arrive at
\[ 
\int  d x_1\left( \sum_\tau \int W_\xi\left[ \left(\begin{array}{cc} w & 0 \\ 0 & I_2 \end{array} \right)
\left(\begin{array}{cc} \tau & 0 \\ 0 & I_2 \end{array} \right)  x_1 x_2  
\right]  d x_2\right) \,.\] 
Here the integral is for $x_1\in X_1(F)\bs X_1(\mathbb{A})$ and $x_2\in X_2(\mathbb{A})$.
We will show that for every $\tau$ and every $\xi$ the following integral is $0$: 
\[ \int _{ X_1(F)\bs X_1(\mathbb{A})} d x_1  \left(\int_{X_2(\mathbb{A})} W_\xi\left[ \left(\begin{array}{cc} w & 0 \\ 0 & I_2 \end{array} \right)
\left(\begin{array}{cc} \tau & 0 \\ 0 & I_2 \end{array} \right)  x_1 x_2  
\right]  d x_2 \right)\,.\]
In order for this expression to even make sense we better show first that, on $X(\mathbb{A})$, the function
\[ x \mapsto \int_{X_2(\mathbb{A})} W_\xi\left[ \left(\begin{array}{cc} w & 0 \\ 0 & I_2 \end{array} \right)
\left(\begin{array}{cc} \tau & 0 \\ 0 & I_2 \end{array} \right)  x x_2  
\right]  d x_2 \]
is invariant under $X_1(F)$. Recall it is invariant under $X_2(\mathbb{A})$. 
So it amounts to the same to prove it is invariant under $X_1(F)X_2(\mathbb{A})$.
Since $\tau$ normalizes the groups $X(\mathbb{A})$ and $X_1(F)X_2(\mathbb{A})$ it amounts to the same to prove that on $X(\mathbb{A})$ the function
\[ x \mapsto \int W_\xi\left[ \left(\begin{array}{cc} w & 0 \\ 0 & I_2 \end{array} \right) x
\left(\begin{array}{cc} \tau & 0 \\ 0 & I_2 \end{array} \right)   x_2  
\right]  d x_2 \]
is invariant under $X_1(F)X_2(\mathbb{A})$. The invariance under $X_2(\mathbb{A})$ being clear we check the invariance under $X_1(F)$. But if $x_1\in X_1(F)$ we have
\[ \int W_\xi\left[ \left(\begin{array}{cc} w & 0 \\ 0 & I_2 \end{array} \right) x_1 x
\left(\begin{array}{cc} \tau & 0 \\ 0 & I_2 \end{array} \right)   x_2  
\right]  d x_2 \] 
\[ = \int W_\xi\left[y_1 \left(\begin{array}{cc} w & 0 \\ 0 & I_2 \end{array} \right) 
 x
\left(\begin{array}{cc} \tau & 0 \\ 0 & I_2 \end{array} \right)   x_2  
\right]  d x_2\]
with
\[ y_1 = \left(\begin{array}{cc} w & 0 \\ 0 & I_2 \end{array} \right) x_1 \left(\begin{array}{cc} w & 0 \\ 0 & I_2 \end{array} \right)^{-1} \,.\]
Since $y_1$ is in $U_{n-2}(F)$, this expression does not depend on $x_1$.

At this point we can reformulate our goal as follows: we have to prove that for every $\tau$ and every $\xi$, the integral of the function
\[ x \mapsto \int W_\xi\left[ \left(\begin{array}{cc} w & 0 \\ 0 & I_2 \end{array} \right)
\left(\begin{array}{cc} \tau & 0 \\ 0 & I_2 \end{array} \right)  x x_2  
\right]  d x_2 \]
over the quotient
\[ X_1(F) X_2(\mathbb{A}) \bs X(\BA)\]
is $0$. Conjugation by $\tau$ defines an automorphism of this quotient which preserves the Haar measure. Hence it suffices to prove that the integral  of the function
\[ x \mapsto \int W_\xi\left[ \left(\begin{array}{cc} w & 0 \\ 0 & I_2 \end{array} \right) x
\left(\begin{array}{cc} \tau & 0 \\ 0 & I_2 \end{array} \right)  x_2  
\right]  d x_2 \] 
over the same quotient vanishes.
Equivalently, we want to prove that
\[ \int_{X_1(F)\bs X_1(\mathbb{A})} d x_1\int_{X_2(\mathbb{A})} W_\xi\left[ \left(\begin{array}{cc} w & 0 \\ 0 & I_2 \end{array} \right) x_1
\left(\begin{array}{cc} \tau & 0 \\ 0 & I_2 \end{array} \right)   x_2  
\right]  d x_2 =0\]
for all $\tau$ and all $\xi$. 
At this point we may exchange the order of integration. 
So it will be enough to prove that
\[ \int_{X_1(F)\bs X_1(\mathbb{A})} d x_1W_\xi\left[ \left(\begin{array}{cc} w & 0 \\ 0 & I_2 \end{array} \right) x_1\right] g = 0
\]
for all $\xi$ and $g$. Now let $Y$ be the subgroup defined by
\[ Y(F) = U_{n-2}(F) \cap w X_1(F) w ^{-1} \,.\]
The integral can be written as
\[ \int _{Y(F)\backslash Y(\mathbb{A})} \psi_{U_n}(y) d y\,\,  W_\xi\left[ \left(\begin{array}{cc} w & 0 \\ 0 & I_2 \end{array} \right)g\right] 
\,.\]
Now we claim that the subgroup $Y$ contains a root subgroup for a positive simple root.
Thus the character $\psi_{U_n}$ is non-trivial on the subgroup $Y(\mathbb{A})$ and the integral vanishes which concludes the proof.

It remains to prove the claim. Since the proof requires an inductive argument we state the claim as a separate lemma.  \qed

\begin{lemma} \label{VanishingLemma}
Let $X$ be the unipotent radical of the standard parabolic subgroup of type $(1,m-1)$ in $GL_m$. Let $w\in GL_m$ be a permutation matrix whose last row has the form $(0,\bullet,\bullet,\ldots,\bullet)$. Then conjugation by $w$ changes one of the root subgroups of $X$ into the root subgroup associated to a positive simple root.
\end{lemma}

\textsc{Proof:} Our assertion is trivial if $m=2$ because then $w=I_2$. So we may assume $m>2$ and our assertion true for $m-1$.
The matrix $w$ has the form
\[ w = \left(\begin{array}{cc} w_1 & 0 \\0&1 \end{array}\right)\left(\begin{array}{cc} 1 & 0 \\ 0 & w_2\end{array}\right) \]
where $w_1$ and $w_2$ are permutation matrices of size $m-1$. Since 
\[ \left(\begin{array}{cc} 1 & 0 \\ 0 & w_2\end{array}\right) \]
normalizes $X$, it suffices to prove our assertion for
\[ w = \left(\begin{array}{cc} w_1 & 0 \\0&1 \end{array}\right) \,.\]
If the last row of $w_1$ has the form
\[ (1,0,\ldots ,0)\,, \]
the root group corresponding to $e_1- e_m$ is conjugated to the root group corresponding to $e_{m-1}-e_m$. So our assertion is true in this case. If the last row of $w_1$ has the form
\[ (0,\bullet,\ldots, \bullet)\,, \]
we can apply the induction hypothesis to $w_1$ and obtain again our assertion. \qed

\textsc{Proof of Theorem \ref{MainResult}:} Recall from \eqref{equvan} that
\[ \int_{(F\backslash \mathbb{A})^3} (U_\xi - V_\xi) \left(\begin{array}{ccc} 1 &  x e_0 & z \\ 0 & I_{n-2} & y f_0 \\ 0 &0 & 1 \end{array} \right) \psi( \alpha x + \beta y)d z d x d y =0 \]
for all $\xi \in V_{\pi}$, $\alpha \in F^\times, \beta \in F^\times$. Here $e_0$ and $f_0$ are of size $n-2$ and
\[ e_0 = (0,0,\ldots,0,1)\,,\, f_0 = \left(\begin{array}{c} 0 \\ 0\\ \vdots \\0 \\ 1 \end{array}\right) \,.\]
Note that if $\Phi_{\xi}=0$ for all $\xi \in V_{\pi}$, then the condition $\CE(n,n-3)$ holds, since the integral over $z$ is just an inner integral of the above integral. Hence, in the following, we assume that the condition $\CE(n,n-3)$ holds, that is, the above integral equals to $0$. 

We can conjugate by a matrix of the form
\[  \left(\begin{array}{ccc} 1 &  0 & 0 \\ 0 & \gamma & 0 \\ 0 &0 & 1 \end{array} \right)\,,\, \gamma \in G_{n-2}(F) \]
to obtain 
\[  \int_{(F\backslash \mathbb{A})^3} (U_\xi - V_\xi) \left(\begin{array}{ccc} 1 &  x e_0 \gamma ^{-1} & z \\ 0 & I_{n-2} &  y \gamma f_0 \\ 0 &0 & 1 \end{array} \right) \psi( \alpha x + \beta y)d z d x d y =0 \,,\]
for all $\xi \in V_{\pi}$, $\alpha \in F^\times, \beta \in F^\times$. 
Thus
\[  \int_{(F\backslash \mathbb{A})^3} (U_\xi - V_\xi) \left(\begin{array}{ccc} 1 &  x e & z \\ 0 & I_{n-2} &  y  f \\ 0 &0 & 1 \end{array} \right) \psi( \alpha x + \beta y)d z d x d y =0 \,,\]
for $e$ (resp. $f$) an $F$-row (resp. column) of size $n-2$ and $\langle e, f \rangle=1$ and for all $\xi \in V_{\pi}$, $\alpha \in F^\times, \beta \in F^\times$. Moreover, right translating by an adelic matrix in the unipotent radical of the parabolic subgroup of type $(1,n-2,1)$ we obtain
\[ \int_{(F\backslash \mathbb{A})^3} \Phi_{\xi} ( u + x e , v+ y f)  \psi( \alpha x + \beta y)d z d x d y =0 \,,\]
for all $(u,v) \in \BA^{n-2} \times \BA^{n-2}$, and for all $\xi \in V_{\pi}$, $\alpha \in F^{\times}$, $\beta \in F^{\times}$. Thus by Proposition \ref{Separable}, 
the function $\Phi_{\xi}$ is separable, for all $\xi \in V_{\pi}$. By Propositions \ref{Vanishing} and \ref{Vanishing2}, it is in fact $0$ and we are done. 
\qed


\begin{thebibliography}{}

\bibitem[Cog02]{Cog02}
J. Cogdell,
{\it L-functions and Converse Theorems for $GL_n$},
Automorphic forms and applications,
97--177, IAS/Park City Math. Ser., \textbf{12}, Amer. Math. Soc., Providence, RI, 2007.

\bibitem[CPS94]{CPS94}
J. Cogdell and I. Piatetski-Shapiro,
{\it Converse theorems for {${GL}_n$}},  Inst. Hautes \'Etudes Sci. Publ. Math. No. \textbf{79} (1994), 157--214.
  
\bibitem[CPS96]{CPS96}
J. Cogdell and I. Piatetski-Shapiro,
{\it A Converse Theorem for $GL_4$},
Math. Res. Letters \textbf{3} (1996), 67--76.
  
\bibitem[CPS99]{CPS99}
J. Cogdell and I. Piatetski-Shapiro,
{\it Converse theorems for {${
  GL}_n$}, {II}}, J. Reine Angew. Math. \textbf{507} (1999), 165--188. 

\bibitem[JPSS79] {JPSS79}
H. Jacquet, I. Piatetski-Shapiro and J. Shalika,
{\it Automorphic forms on $GL(3)$.}
Ann. of Math. (2) \textbf{109} (1979), no. 1--2, 169--258.

\bibitem[JL70] {JL70}
H. Jacquet and R. Langlands,
{\it Automorphic forms on $GL(2)$}. Lecture Notes in Mathematics, Vol. \textbf{114}. Springer-Verlag, Berlin-New York, 1970. 

\bibitem[PS76] {PS76}
I. Piatetski-Shapiro,
{\it Zeta functions of $GL(n)$}. 
Univ. Maryland Techn. Rep. TR 76--46 (1976).
\end{thebibliography}
\end{document}